\documentclass{amsart}
\usepackage{amsfonts}

\usepackage{graphicx}

\input{tcilatex}
\input tcilatex

\begin{document}
\title[Sibuya progenies]{On Sibuya trees and forests}
\author{Thierry E. Huillet}
\address{Laboratoire de Physique Th\'{e}orique et Mod\'{e}lisation \\
CY Cergy Paris University, CNRS UMR-8089 \\
2 Avenue Adolphe Chauvin, 95302, Cergy-Pontoise, FRANCE\\
E-mail: Thierry.Huillet@cyu.fr}
\maketitle

\begin{abstract}
We show that the Sibuya distribution and its non-critical relatives are
relevant in the context of the recursive generation of both simply generated
and increasing critical trees' and forests' progenies. A special class of
generalized Stirling numbers are at the heart of the analysis of the induced
occupancy distributions. Asymptotic aspects of large forests are addressed.%
\newline



Journal of the Indian Society for Probability and Statistics.
\end{abstract}

\section{Introduction}

The discrete Sibuya distribution on the non-negative integers introduced in
Sibuya (\textsl{1979}) is a particular version of the generalized Sibuya
distribution appearing and studied in Dacey (\textsl{1969}); Kozubowski and
Podg\'{o}rski (\textsl{2018}) and Huillet (\textsl{2020}). The latter is a
class of distributions with Gauss hypergeometric probability generating
functions (pgf's). It includes Yule (\textsl{1925}) and Simon (\textsl{1955}%
) distributions.

In the context of Bienaym\'{e}-Galton-Watson (BGW) processes, it appears as
the law of its offspring branching mechanism [Huillet (\textsl{2016}), for
example], one important feature in this context being the stability under
iteration of its pgf $\phi \left( z\right) =1-\lambda \left( 1-z\right)
^{\alpha }$, $\alpha \in \left( 0,1\right) $, $\lambda \in \left( 0,1\right]
.$

Recently, Letac (\textsl{2019}) showed that the Sibuya distribution also
appears as a progeny in a critical BGW process when $\alpha \in \left[
1/2,1\right) $. While considering rescalings, sub-critical and
super-critical versions of the critical Sibuya distribution are herewith
introduced and studied.

We will next show that, for increasing (else recursive) trees' progenies,
appearing in phylogenetic contexts, the Sibuya distribution makes sense for
the full admissible range of $\alpha $.

We will also consider $n-$forests of Sibuya such trees, so both simply
generated and increasing, emphasizing their recursive generation while $%
n\rightarrow n+1$. We will study the number of trees entering in the
composition of a size$-n$ forest of such trees, the joint sizes of the
constitutive trees in a size-$n$ forest of $k$ trees and the one-dimensional
marginal size of a typical tree. We derive asymptotic results when $%
n\rightarrow \infty $, $k\rightarrow \infty ,$ separately, or when $n$,$%
k\rightarrow \infty $, jointly, while $n/k\rightarrow \rho >1,$ in a
`thermodynamic' limit.

Both approaches rest on the analysis of generating functions and can
sometimes take advantage of the Lagrange inversion formula.

\section{Sibuya distribution as a BGW progeny}

With $\alpha \in \left( 0,1\right) $ and $z\in \left[ 0,z_{c}\right] $ where 
\begin{equation*}
z_{c}:=\sup \left\{ z>0:\Phi \left( z\right) <\infty \right\} =1,
\end{equation*}
consider the pgf of a Sibuya random variable (rv) on the positive integers 
\begin{equation*}
\Phi \left( z\right) =1-\left( 1-z\right) ^{\alpha }=\alpha z\cdot F\left(
1,1-\alpha ;2;z\right) ,
\end{equation*}
where $F={_{2}}F_{1}$ is the Gauss hypergeometric function. Letac (\textsl{%
2019}) shows that it is the one of the total progeny of a rooted BGW process
with one founder when $\alpha \in \left[ 1/2,1\right) $ (and shows that it
is not a progeny if $\alpha \in \left( 0,1/2-10^{-9}\right) $ while
conjecturing that this is also true when $\alpha \in \left( 0,1/2\right) $).
Indeed, in that case, see Harris (\textsl{1963}), $\Phi \left( z\right) $
solves the functional equation $\Phi \left( z\right) =z\phi \left( \Phi
\left( z\right) \right) $, $\Phi \left( 0\right) =0$ where $\phi \left(
z\right) $ is the regular ($\phi \left( 1\right) =1$) branching mechanism 
\begin{equation}
\phi \left( z\right) =\frac{z}{1-\left( 1-z\right) ^{1/\alpha }},
\label{FF0}
\end{equation}
with $\left[ z^{n}\right] \phi \left( z\right) >0$ and $\phi \left( 0\right)
=\alpha .$ Here $\left[ z^{n}\right] \phi \left( z\right) $ the $z^{n}-$%
coefficient in the power-series expansion of $\phi \left( z\right) $. Note 
\begin{equation*}
z_{*}:=\sup \left\{ z>0:\phi \left( z\right) <\infty \right\} =z_{c}=1,
\end{equation*}
so that both convergence radii of $\phi \left( z\right) $ and $\Phi \left(
z\right) $ coincide. When $\alpha =1/2$, $\phi \left( z\right) =1/\left(
2-z\right) ,$ the pgf of a shifted geometric$\left( 1/2\right) $
distribution.

The underlying BGW process is critical with extinction probability $1$,
observing $\phi ^{\prime }\left( 1\right) =1$. Denoting $N\left( 1\right) $
the one-founder total progeny of this BGW process, then 
\begin{equation*}
\mathbf{P}\left( N\left( 1\right) =n\right) =\left[ z^{n}\right] \Phi \left(
z\right) =\frac{\alpha \left[ 1-\alpha \right] _{n-1}}{n!}\underset{%
n\rightarrow \infty }{\sim }\alpha n^{-\left( \alpha +1\right) }/\Gamma
\left( 1-\alpha \right)
\end{equation*}
(a power-law with index $\beta :=\alpha +1\in \left( 1,2\right) $). Here $%
\left[ a\right] _{b}=a\left( a+1\right) \cdots \left( a+b-1\right) ,$ as the
rising factorial of $a$. In particular, we have $\Phi ^{\prime }\left(
1\right) =\mathbf{E}\left( N\left( 1\right) \right) =\infty $.\newline

\textsl{Remark:} Lagrange inversion formula [Stanley (\textsl{1999}) and
Surya and Warnke (\textsl{2023}), for example] states that for all $n\geq 1,$%
\begin{equation*}
\left[ z^{n}\right] \Phi \left( z\right) =\frac{1}{n}\left[ z^{n-1}\right]
\phi \left( z\right) ^{n}.
\end{equation*}
Hence here the non-obvious identity, 
\begin{equation*}
\frac{1}{n}\left[ z^{n-1}\right] \left[ \frac{z}{1-\left( 1-z\right)
^{1/\alpha }}\right] ^{n}=\frac{\alpha \left[ 1-\alpha \right] _{n-1}}{n!}.
\end{equation*}

A more general form of Lagrange inversion formula states that (with $%
^{\prime }$ denoting derivative) 
\begin{equation}
\left[ z^{n}\right] h\left( \Phi \left( z\right) \right) =\frac{1}{n}\left[
z^{n-1}\right] \left( h^{\prime }\left( z\right) \phi \left( z\right)
^{n}\right) .  \label{F3a}
\end{equation}
for any arbitrary analytic output function $h$.

- If $h\left( z\right) $ is a pgf, $h\left( \Phi \left( z\right) \right) $
is the pgf of the total progeny of a branching process generated by $\phi $
with a random number of founders, say $N_{0},$ for which $h\left( z\right) =%
\mathbf{E}\left( z^{N_{0}}\right) .$

- If $h\left( z\right) =\left( 1-uz\right) ^{-1}$ where $u$ `marks' the
number of distinguishable trees in a forest of simple random trees generated
by $\phi $, with $K\left( z,u\right) :=1/\left( 1-u\Phi \left( z\right)
\right) $%
\begin{equation}
\left[ z^{n}\right] K\left( z,u\right) =\frac{u}{n}\left[ z^{n-1}\right]
\left( \left( 1-uz\right) ^{-2}\phi \left( z\right) ^{n}\right) .
\label{F3b}
\end{equation}
This consists of a degree-$n$ polynomial in $u$. The double generating
function $K\left( z,u\right) $ may be viewed as a `grand-canonical'
partition function. In that case, 
\begin{equation*}
\frac{\left[ u^{k}z^{n}\right] K\left( z,u\right) }{\left[ z^{n}\right]
K\left( z,1\right) }=\frac{\frac{k}{n}\left[ z^{n-k}\right] \phi \left(
z\right) ^{n}}{\frac{1}{n}\left[ z^{n-1}\right] \left( \left( 1-z\right)
^{-2}\phi \left( z\right) ^{n}\right) }=\mathbf{P}\left( K_{n}=k\right) ,
\end{equation*}
where $K_{n}$ is the number of distinguishable (labelled) trees forming the
size$-n$ forest. $K\left( z,u\right) $ consist of a `sequence' of trees
forming the forest, with $u$ `marking' their number.

- Would the trees be indistinguishable (we will not consider this case), the
same formulae hold but now with $E\left( z,u\right) :=e^{u\Phi \left(
z\right) },$ consisting of a `set' of trees forming the forest, so with 
\begin{equation}
\left[ u^{k}z^{n}\right] E\left( z,u\right) =\frac{1}{k!}\frac{k}{n}\left[
z^{n-k}\right] \phi \left( z\right) ^{n}.  \label{F3c}
\end{equation}
The use of the 'set' partition function $E\left( z,u\right) $, instead of
the 'sequenced' partition function $K\left( z,u\right) $, applies when
considering the number of partitions of $n$ distinguishable atoms into $k$
indistinguishable parts.

When conditioning on the size of the forest of distinguishable trees, the
joint Maxwell-Boltzmann exchangeable occupancy law of a population of $k$
clusters in a size-$n$ forest, as a partition structure, is 
\begin{equation}
\left\{ 
\begin{array}{c}
\mathbf{P}\left( \mathbf{N}_{n,k}=\mathbf{n}_{k}\right) =\frac{\left[
z^{n}\prod_{l=1}^{k}z_{l}^{n_{l}}\right] \prod_{l=1}^{k}\Phi \left(
zz_{l}\right) }{\left[ z^{n}\right] \Phi \left( z\right) ^{k}} \\ 
=\frac{\prod_{l=1}^{k}\mathbf{P}\left( N\left( 1\right) =n_{l}\right) }{%
\mathbf{P}\left( N\left( k\right) =n\right) }\delta _{\left| \mathbf{n}%
_{k}\right| =n} \\ 
=\frac{n!}{\mathcal{C}_{n,k}}.\alpha ^{k}\prod_{l=1}^{k}\frac{\left[
1-\alpha \right] _{n_{l}-1}}{n_{l}!}\delta _{\left| \mathbf{n}_{k}\right|
=n}.
\end{array}
\right.  \label{two}
\end{equation}
Here $N\left( k\right) =\sum_{l=1}^{k}N\left( 1\right) _{l}$, where the $%
N\left( 1\right) _{l}$ are i.i.d. copies of $N\left( 1\right) $, $\mathcal{C}%
_{n,k}=n!\left[ z^{n}\right] \Phi \left( z\right) ^{k}$ the normalization
factor and $\mathbf{n}_{k}=\left( n_{1},\cdots ,n_{k}\right) $ positive
integers summing to $n$, say $\left| \mathbf{n}_{k}\right| =n$, indicating
the values of the random box occupancies $\mathbf{N}_{n,k}=\left(
N_{n,1}\left( k\right) ,\cdots ,N_{n,k}\left( k\right) \right) $.

\textbf{The distribution of the size of a typical (1-dimensional marginal)
box occupancy} is obtained by summation upon the other variables, and it
reads 
\begin{equation}
\mathbf{P}\left( N_{n,1}\left( k\right) =n_{1}\right) =\frac{\left[
z^{n_{1}}\right] \Phi \left( z\right) \left[ z^{n-n_{1}}\right] \Phi \left(
z\right) ^{k-1}}{\left[ z^{n}\right] \Phi \left( z\right) ^{k}}
\label{DDSmarg}
\end{equation}

\begin{equation*}
=\binom{n}{n_{1}}\frac{\mathcal{C}_{n_{1},1}\cdot \mathcal{C}_{n-n_{1},k-1}}{%
\mathcal{C}_{n,k}}.
\end{equation*}
Here, $\mathcal{C}_{n_{1},1}=\alpha \left[ 1-\alpha \right] ^{n_{1}-1}$ and
the $\mathcal{C}_{n,k}$'s are given below. With $n_{1}=1,\cdots ,n-k+1$,
summing over $n_{1}$, the probability mass function (pmf) of $N_{n,1}\left(
k\right) $ is also seen non-defective and proper. Clearly, $\mathbf{E}\left(
N_{n,1}\left( k\right) \right) =n/k.$

A \textbf{large population} thermodynamic limit exists ($n,k\rightarrow
\infty $ with $n/k\rightarrow \rho \geq 1$), with $N_{1}=:N_{1}\left( \rho
\right) $ having the mean-$\rho $ limiting distribution given, from a
saddle-point analysis [see for example Bialas et al. (\textsl{2023}),
section III, for a detailed justification], by the regular pgf 
\begin{equation}
\Phi _{\rho }\left( z\right) :=\mathbf{E}\left( z^{N_{1}\left( \rho \right)
}\right) =\frac{\Phi \left( zz_{\rho }\right) }{\Phi \left( z_{\rho }\right) 
}=\frac{1-\left( 1-zz_{\rho }\right) ^{\alpha }}{1-\left( 1-z_{\rho }\right)
^{\alpha }},  \label{f19}
\end{equation}
where $z_{\rho }\in \left( 0,z_{c}\right) $ solves implicitly 
\begin{equation}
\Psi \left( z\right) :=\frac{z\Phi ^{\prime }\left( z\right) }{\Phi \left(
z\right) }=\frac{\alpha z\left( 1-z\right) ^{-\left( 1-\alpha \right) }}{%
1-\left( 1-z\right) ^{\alpha }}=\rho .  \label{f20}
\end{equation}
This pgf is the one of a canonical exponential family.

The free energy per tree in the thermodynamic limit is 
\begin{eqnarray*}
f\left( \rho \right) &=&\lim_{k\rightarrow \infty }-\frac{1}{k}\log \left[
z^{n}\right] \Phi \left( z\right) ^{k}\text{, }n/k\rightarrow \rho \\
&=&z_{\rho }\log z_{\rho }-\log \Phi \left( z_{\rho }\right) .
\end{eqnarray*}
For all $\rho \geq 1$, $z_{\rho }$ is uniquely defined firstly because $\Psi
\left( z\right) $ is increasing, so one-to-one, and because, when $%
z\rightarrow z_{c}=1$, $\Phi \left( z\right) \rightarrow \Phi \left(
z_{c}\right) =1$ and $\Phi ^{\prime }\left( z\right) \rightarrow \infty ,$
with $\rho \rightarrow \infty $, so that $\Psi \left( z_{c}\right) =\infty .$
This shows that the Sibuya model does not show up the condensation
phenomenon discussed in Janson (\textsl{2012}) when $\Psi \left(
z_{c}\right) <\infty $, so that $z_{\rho }$ would be defined only for $\rho
\leq \rho _{c}:=\Psi \left( z_{c}\right) .$ A condensation phenomenon (with
a giant component) would indeed occur if $\left[ z^{n}\right] \Phi \left(
z\right) \underset{n\rightarrow \infty }{\sim }C\cdot n^{-\beta }$ with $%
\beta >2$, as is the case, for instance, for the zeta series weight model
[see Janson (\textsl{2012}) and Bialas et al. (\textsl{2023})].

In the thermodynamic limit, the Sibuya forest is made of a countable number
of i.i.d. trees, the size of each of which, say $N_{l}\left( \rho \right) $, 
$l\geq 1$, has distribution given by (\ref{f19}). By Cram\'{e}r's theorem (%
\textsl{1938}), with $\overline{N}_{\rho }\left( k\right)
:=\sum_{l=1}^{k}N_{l}\left( \rho \right) $ and $r>1,$%
\begin{equation}
-\frac{1}{k}\log \mathbf{P}\left( \overline{N}_{\rho }\left( k\right)
/k\rightarrow r\right) \underset{k\rightarrow \infty }{\rightarrow }f_{\rho
}\left( r\right) \geq 0.  \label{f16}
\end{equation}
The rate function $f_{\rho }\left( r\right) $ is the Legendre transform of
the convex function $\log \Phi _{\rho }\left( z\right) $, hence, 
\begin{equation}
f_{\rho }\left( r\right) =r\log z_{r}-\log \Phi _{\rho }\left( z_{r}\right) 
\text{ where }\frac{z_{r}\Phi _{\rho }^{\prime }\left( z_{r}\right) }{\Phi
_{\rho }\left( z_{r}\right) }=\Psi \left( z_{\rho }z_{r}\right) =r.
\label{f2}
\end{equation}
Note $f_{\rho }\left( \rho \right) =0$ as a result of $z_{r}=1$ when $r=\rho 
$, translating $\overline{N}_{\rho }\left( k\right) /k\rightarrow \rho ,$
almost surely.

We will show below that $\Phi _{\rho }\left( z\right) $ is the modified
(rescaled) regular pgf of the progeny of a regular subcritical BGW process.

\subsection{$k$-forests of distinguishable Sibuya trees}

We now observe from (21) in Huillet (\textsl{2024}) that, by a Faa di Bruno
formula, 
\begin{equation}
\mathcal{S}_{n,k}:=S_{n,k}\left( -1,-\alpha ,0\right) =\frac{n!}{k!}%
\sum^{*}\prod_{l=1}^{k}\frac{\left[ 1-\alpha \right] _{n_{l}-1}}{n_{l}!},
\label{faadibruno}
\end{equation}
where the star sum runs over the positive sample sizes $\left(
n_{l};l=1,\cdots ,k\right) $ obeying $\sum_{l=1}^{k}n_{l}=n$ (the number of
compositions of $n$ into $k$ parts; there are $\binom{n-1}{k-1}$ terms in
this sum)$.$ The $\mathcal{S}_{n,k}$ are then particular generalized
Stirling numbers of Hsu and Shiue (\textsl{1998}), in the notations of
Huillet (\textsl{2024}). As a result, from (\ref{two}): $\mathcal{C}%
_{n,k}=k!\alpha ^{k}\mathcal{S}_{n,k}.$ See Appendix \textbf{A1} for a recap.%
\newline

\textsl{Remark:} Lagrange inversion formula states that for all $n\geq 1$%
\begin{equation}
\left[ z^{n}\right] \Phi \left( z\right) ^{k}=\frac{k}{n}\left[
z^{n-k}\right] \phi \left( z\right) ^{n}.  \label{FF3}
\end{equation}
Hence here from (\ref{faadibruno}), for a forest of distinguishable Sibuya
trees: 
\begin{equation}
\frac{k}{n}\left[ z^{n-k}\right] \left[ \frac{z}{1-\left( 1-z\right)
^{1/\alpha }}\right] ^{n}=\frac{\mathcal{C}_{n,k}}{n!}:=\frac{k!}{n!}\alpha
^{k}\mathcal{S}_{n,k}.  \label{FF4}
\end{equation}

The $\mathcal{S}_{n,k}$ obey the three-terms recurrence (Stirling's
triangle) 
\begin{equation}
\mathcal{S}_{n+1,k}=\mathcal{S}_{n,k-1}+\left( n-k\alpha \right) \mathcal{S}%
_{n,k},  \label{rec}
\end{equation}
with boundary conditions $\mathcal{S}_{n,0}=\delta _{n,0}$ and $\mathcal{S}%
_{n,n}=1$ for all $n\geq 0$. The $\mathcal{S}_{n,k}$ are all positive for $%
k=1,\cdots ,n$ and for all $\alpha \in \left( 0,1\right) .$ Furthermore,
with $\left( a\right) _{n}:=a\left( a-1\right) \cdots \left( a-n+1\right) $%
\begin{equation*}
\mathcal{S}_{n,k}=\frac{\alpha ^{-k}}{k!}\sum_{l=0}^{k}\left( -1\right)
^{n+l}\binom{k}{l}\left( l\alpha \right) _{n}
\end{equation*}
is an explicit representation in terms of an alternating sum (see Appendix 
\textbf{A1}).

Hence, $\mathcal{C}_{n,k}=k!\alpha ^{k}\mathcal{S}_{n,k}$ obeys the
recurrence 
\begin{equation*}
\mathcal{C}_{n+1,k}=k\alpha \mathcal{C}_{n,k-1}+\left( n-k\alpha \right) 
\mathcal{C}_{n,k}.
\end{equation*}
With $\mathcal{C}_{n}=\sum_{k=1}^{n}\mathcal{C}_{n,k}=n!\left[ z^{n}\right] 
\frac{1}{1-\Phi \left( z\right) }=n!\left[ z^{n}\right] \left( 1-z\right)
^{-\alpha }=\left[ \alpha \right] _{n},$ the total number of size-$n$
forests with distinguishable Sibuya trees, $\mathbf{P}\left( K_{n}=k\right) =%
\frac{\mathcal{C}_{n,k}}{\mathcal{C}_{n}}$ thus obeys the triangular
recurrence ($k=1,\cdots ,n$) 
\begin{equation}
\mathbf{P}\left( K_{n+1}=k\right) =\frac{k\alpha }{\alpha +n}\mathbf{P}%
\left( K_{n}=k-1\right) +\frac{n-k\alpha }{\alpha +n}\mathbf{P}\left(
K_{n}=k\right) ,  \label{FF5}
\end{equation}
with $K_{0}=0$ and $K_{1}=1$. The sequence $\left( K_{n}\right) _{n}$ is a
pure birth Markov chain on the state-space $\mathbf{N}:=\left\{ 1,2,\cdots
\right\} $. The transition probabilities from state $k$ to states $\left\{
k+1,k\right\} $ do sum to $1$: $\frac{\left( k+1\right) \alpha }{\alpha +n}+%
\frac{n-k\alpha }{\alpha +n}=1$. Low order iterations of (\ref{FF5}) yield
for instance 
\begin{equation*}
\mathbf{P}\left( K_{2}=1\right) =\frac{1-\alpha }{\alpha +1},\text{ }\mathbf{%
P}\left( K_{2}=2\right) =\frac{2\alpha }{\alpha +1},
\end{equation*}
translating a propensity to form two singletons if $\alpha >1/3$ (as a
result of $\mathbf{P}\left( K_{2}=2\right) >\mathbf{P}\left( K_{2}=1\right) $%
) and

\begin{eqnarray*}
\mathbf{P}\left( K_{3}=1\right) &=&\frac{\left( 1-\alpha \right) \left(
2-\alpha \right) }{\left( \alpha +1\right) \left( \alpha +2\right) }, \\
\mathbf{P}\left( K_{3}=2\right) &=&\frac{6\alpha \left( 1-\alpha \right) }{%
\left( \alpha +1\right) \left( \alpha +2\right) }\text{,} \\
\mathbf{P}\left( K_{3}=3\right) &=&\frac{6\alpha ^{2}}{\left( \alpha
+1\right) \left( \alpha +2\right) }.
\end{eqnarray*}

Multiplying (\ref{FF5}) by $k$ and summing over all possible $k\in \left\{
1,\cdots ,n+1\right\} $ shows that the mean $\mu _{n+1}$ of $K_{n+1}$ obeys (%
$\mu _{1}=1$) 
\begin{eqnarray*}
\left( n+\alpha \right) \mu _{n+1} &=&\alpha \left( 2\mu _{n}+1\right) +n\mu
_{n} \\
\left( n+\alpha \right) \left( \mu _{n+1}-\mu _{n}\right) &=&\alpha \left(
1+\mu _{n}\right)
\end{eqnarray*}
Approximating heuristically this difference equation in discrete time by the
differential equation in continuous time ($\mu _{1}=1$): 
\begin{equation*}
\left( t+\alpha \right) \overset{.}{\mu }_{t}=\alpha \left( 1+\mu
_{t}\right) ,
\end{equation*}
we get the solution: $\mu _{t}=2\left( \frac{t+\alpha }{1+\alpha }\right)
^{\alpha }-1$. If indeed the differential equation turns out to be a good
approximation to the difference equation, we would be able to infer that, as 
$n\rightarrow \infty ,$%
\begin{equation*}
\mu _{n}\sim C_{\alpha }\cdot n^{\alpha },
\end{equation*}
for some prefactor $C_{\alpha }>0.$

The rv $N\left( 1\right) $ is heavy-tailed with tail index $\alpha $ and $%
N\left( k\right) $ is the sum of $k$ i.i.d. such rv's. The process $\left(
N\left( k\right) \right) _{k}$ has stationary independent increments and the
processes $\left( N\left( k\right) \right) _{k}$ and $\left( K_{n}\right)
_{n}$ are mutual inverses in that 
\begin{equation*}
K_{n}=\inf \left( k\geq 1:\text{ }N\left( k\right) >n\right) .
\end{equation*}
The process $\left( K_{n}\right) _{n}$ is then a renewal process with times
elapsed between consecutive moves up by one unit all distributed like $%
N\left( 1\right) :$%
\begin{equation*}
K_{n}\overset{d}{=}1\cdot 1_{\left\{ \overline{N}\left( 1\right)
>n-1\right\} }+\sum_{m=1}^{n-1}1_{\left\{ N\left( 1\right) =m\right\} }\cdot
\left( 1+K_{n-m}^{\prime }\right) \text{, }n\geq 1.
\end{equation*}

With $M_{n}:=\frac{\Gamma \left( n+\alpha \right) }{\Gamma \left( n+2\alpha
\right) }\left( K_{n}+1\right) ,$ it can be checked from (\ref{FF5}) that
the process $\left( M_{n}\right) _{n}$ is a non-negative martingale
converging almost surely to a positive rv, say $W_{\alpha ,\alpha }$.
Observing $M_{n}\underset{n\rightarrow \infty }{\sim }n^{-\alpha }K_{n}$, we
get 
\begin{equation}
n^{-\alpha }K_{n}\underset{n\rightarrow \infty }{\rightarrow }W_{\alpha
,\alpha }.  \label{FF6}
\end{equation}
Based on Uchaikin-Zolotarev (\textsl{1999}) p. 62, we have $N\left( k\right)
/k^{1/\alpha }\overset{d}{\rightarrow }S_{\alpha }$\ (a one-sided $\alpha -$%
stable rv): as $k$\ is large indeed, for $\lambda \geq 0,$%
\begin{equation*}
\mathbf{E}e^{-\lambda k^{-1/\alpha }N\left( k\right) }=\left( 1-\left(
1-e^{-\lambda k^{-1/\alpha }}\right) ^{\alpha }\right) ^{k}\sim \left( 1-%
\frac{\lambda ^{\alpha }}{k}\right) ^{k}\sim e^{-\lambda ^{\alpha }}.
\end{equation*}
As a result, $K_{n}/n^{\alpha }\overset{d}{\rightarrow }W_{\alpha ,\alpha }$
the size-biased version of the standard Mittag-Leffler rv $W_{\alpha
,1}:=S_{\alpha }^{-\alpha }.$\emph{\ }The rv $W_{\alpha ,\alpha }$ has a
one-parameter Mittag-Leffler distribution with moment function (see Appendix 
\textbf{A2}) 
\begin{equation*}
\mathbf{E}\left( W_{\alpha ,\alpha }^{q}\right) =\frac{\Gamma \left( \alpha
\right) \Gamma \left( q+1\right) }{\Gamma \left( q\alpha +\alpha \right) },%
\text{ }q>-1,
\end{equation*}
having mean $\frac{\Gamma \left( \alpha \right) }{\Gamma \left( 2\alpha
\right) }$ and variance $\frac{2\Gamma \left( \alpha \right) }{\Gamma \left(
3\alpha \right) }-\left( \frac{\Gamma \left( \alpha \right) }{\Gamma \left(
2\alpha \right) }\right) ^{2}$.

\textsl{Remark}: Although arising from different processes, (\ref{two})
shows some analogy with the occupancy EPPF Pitman (\textsl{2006}), (3.6) p.
63 of the Ewens-Pitman sampling formula from the two-parameters
Poisson-Dirichlet partition $\mathrm{PD}\left( \alpha ,\theta \right) $,
appearing in the Chinese Restaurant Problem. In our one-parameter branching
process approach with $\alpha \in \left[ 1/2,1\right) $, given $n$ species
split into $k$ types or genera, a new incoming species will become the root
species of some new genus with probability $\frac{\left( k+1\right) \alpha }{%
\alpha +n}.$ See the Appendix \textbf{A2} for a detailed connection between
the two occupancy distributions.

\subsection{Rescalings}

Following Meir and Moon (\textsl{1978}) ideas, we herewith consider weighted
versions of Sibuya trees, thereby unveiling the sub-critical and
super-critical counterparts to the critical Sibuya trees just defined.

With $c_{1},c_{2}>0$, consider the weight sequence $w_{b}=c_{1}^{b}c_{2}$
assigned to each node with out-degree $b$ in a size-$n$ tree $\tau _{n}$.
Then, with $n_{b}\left( \tau _{n}\right) $ the number of these nodes in $%
\tau _{n}$, $\prod_{b\geq 0}w_{b}^{n_{b}\left( \tau _{n}\right)
}=c_{1}^{n-1}c_{2}^{n}$, in view of the relations 
\begin{equation*}
\left\{ 
\begin{array}{c}
\sum_{b\geq 0}n_{b}\left( \tau _{n}\right) =n \\ 
\sum_{b\geq 1}bn_{b}=n-1.
\end{array}
\right.
\end{equation*}
As a result, 
\begin{equation*}
\mathbf{P}\left( N\left( 1\right) =n\right) \rightarrow \widetilde{\mathbf{P}%
}\left( N\left( 1\right) =n\right) =c_{1}^{-1}\mathbf{P}\left( N\left(
1\right) =n\right) \left( c_{1}c_{2}\right) ^{n}
\end{equation*}
is the weighted version of $\mathbf{P}\left( N\left( 1\right) =n\right) $.
Equivalently 
\begin{equation}
\Phi \left( z\right) \rightarrow \widetilde{\Phi }\left( z\right)
=c_{1}^{-1}\Phi \left( c_{1}c_{2}z\right) ,  \label{scale}
\end{equation}
solving $\widetilde{\Phi }\left( z\right) =z\widetilde{\phi }\left( 
\widetilde{\Phi }\left( z\right) \right) $ where $\widetilde{\phi }\left(
z\right) =c_{2}\phi \left( c_{1}z\right) $ is the modified `branching
mechanism', not necessarily a regular pgf. It is a regular pgf when $%
c_{2}=1/\phi \left( c_{1}\right) $. Note that, in the case of Sibuya model, $%
\alpha $ still needs to belong to $\left[ 1/2,1\right) $ for $\widetilde{%
\mathbf{P}}\left( N\left( 1\right) =n\right) =\left[ z^{n}\right] \widetilde{%
\Phi }\left( z\right) >0.$

If $c_{1}=1$ and $c_{2}\leq z_{c}$, $\widetilde{\Phi }\left( z\right) =\Phi
\left( c_{2}z\right) $ resulting in a weighted version of $\Phi \left(
z\right) $ with shifted convergence radius $z_{c}\rightarrow \widetilde{z}%
_{c}=z_{c}/c_{2}\geq 1$ [$=1$ if in addition $c_{2}=z_{c}$]$.$ Note that $%
\widetilde{\phi }\left( z\right) =c_{2}\phi \left( z\right) $ no longer is a
branching mechanism if $\phi \left( z\right) $ is one with $\phi \left(
1\right) =1$.

If $c_{1}c_{2}=1$, $\widetilde{\phi }\left( z\right) =c_{1}^{-1}\phi \left(
c_{1}z\right) $ is the modified `branching mechanism', not necessarily a pgf
unless $c_{1}=\phi \left( c_{1}\right) $. Then, $\widetilde{\Phi }\left(
z\right) =c_{1}^{-1}\Phi \left( z\right) $ resulting in a scaled version of $%
\Phi \left( z\right) $ with unmodified convergence radius. If $\Phi \left(
1\right) <1$, choosing $c_{1}=c_{2}^{-1}=\Phi \left( 1\right) $ yields $%
\widetilde{\Phi }\left( 1\right) =1.$ We have that $\widetilde{\Phi }\left(
1\right) =1$ entails $\widetilde{\phi }\left( 1\right) =1$ and, if in
addition, $\widetilde{\Phi }^{\prime }\left( 1\right) <\infty $ then $%
\widetilde{\phi }^{\prime }\left( 1\right) <1$ (a sub-critical case), with $%
\widetilde{\Phi }^{\prime }\left( 1\right) =1/\left( 1-\widetilde{\phi }%
^{\prime }\left( 1\right) \right) $.

We now draw some additional conclusions in case of the specific Sibuya
progeny:

$\bullet $ Scaled Sibuya $c_{1}=1,$ $c_{2}\in \left( 0,1\right) .$ Then $%
\widetilde{\Phi }\left( z\right) =\Phi \left( c_{2}z\right) =1-\left(
1-c_{2}z\right) ^{\alpha }$ with 
\begin{equation*}
\widetilde{\phi }\left( z\right) =c_{2}\phi \left( z\right) =\frac{c_{2}z}{%
1-\left( 1-z\right) ^{1/\alpha }}
\end{equation*}
and $\widetilde{\phi }\left( 1\right) =$ $\widetilde{\phi }^{\prime }\left(
1\right) =c_{2}<1$. Note $\widetilde{\phi }\left( 0\right) =c_{2}\alpha <1$
is the modified probability of no offspring. So, $\widetilde{\Phi }\left(
z\right) $ is the pgf of a super-critical BGW process with extinction
probability $\rho _{e}=\widetilde{\Phi }\left( 1\right) =1-\left(
1-c_{2}\right) ^{\alpha }<1$, which is also the smallest solution in $\left(
0,1\right) $ to $\widetilde{\phi }\left( \rho _{e}\right) =\rho _{e}.$
Because $\widetilde{\phi }\left( 1\right) <1$, the branching mechanism is
not regular having a probability $1-c_{2}$ to take the value $\infty $,
resulting in the mean value of the offspring distribution being infinite.

$\bullet $ Scaled Sibuya $c_{2}=1,$ $c_{1}\in \left( 0,1\right) .$

Then $\widetilde{\Phi }\left( z\right) =c_{1}^{-1}\Phi \left( c_{1}z\right)
=c_{1}^{-1}\left[ 1-\left( 1-c_{1}z\right) ^{\alpha }\right] $ with 
\begin{equation*}
\widetilde{\phi }\left( z\right) =\phi \left( c_{1}z\right) =\frac{c_{1}z}{%
1-\left( 1-c_{1}z\right) ^{1/\alpha }}
\end{equation*}
and $\widetilde{\phi }\left( 1\right) =$ $c_{1}/\left[ 1-\left(
1-c_{1}\right) ^{1/\alpha }\right] <1$. Note $\widetilde{\phi }\left(
0\right) =\alpha $. So, $\widetilde{\Phi }\left( z\right) $ is the pgf of a
super-critical BGW process with extinction probability $\rho _{e}=\widetilde{%
\Phi }\left( 1\right) =c_{1}^{-1}\left[ 1-\left( 1-c_{1}\right) ^{\alpha
}\right] <1$, which is also the smallest solution in $\left( 0,1\right) $ to 
$\widetilde{\phi }\left( \rho _{e}\right) =\rho _{e}.$ Because $\widetilde{%
\phi }\left( 1\right) <1$, the modified branching mechanism is not regular
having a probability $1-\rho _{e}$ to take the value $\infty $: the mean
value of the offspring distribution is infinite.

$\bullet $ Scaled Sibuya $c_{1}\in \left( 0,1\right) $, $c_{2}>0$, $%
c_{1}c_{2}<1.$

Then $\widetilde{\Phi }\left( z\right) =c_{1}^{-1}\Phi \left(
c_{1}c_{2}z\right) =c_{1}^{-1}\left[ 1-\left( 1-c_{1}c_{2}z\right) ^{\alpha
}\right] $ with 
\begin{equation*}
\widetilde{\phi }\left( z\right) =c_{2}\phi \left( c_{1}z\right) =\frac{%
c_{1}c_{2}z}{1-\left( 1-c_{1}z\right) ^{1/\alpha }}.
\end{equation*}
$\widetilde{\phi }\left( z\right) $ is a regular pgf ($\widetilde{\phi }%
\left( 1\right) =1$) iff $c_{2}=c_{1}^{-1}\left[ 1-\left( 1-c_{1}\right)
^{1/\alpha }\right] >1$. In that case, $\widetilde{\Phi }\left( 1\right) =1$
and extinction is almost sure. In addition, $\widetilde{\phi }^{\prime
}\left( 1\right) =1-\frac{c_{1}}{\alpha }\frac{\left( 1-c_{1}\right)
^{1/\alpha -1}}{1-\left( 1-c_{1}\right) ^{1/\alpha }}<1:$ the model is
sub-critical.

We have $\rho _{e}=\widetilde{\Phi }\left( 1\right) =c_{1}^{-1}\left[
1-\left( 1-c_{1}c_{2}\right) ^{\alpha }\right] <1$ iff $c_{2}<c_{1}^{-1}%
\left[ 1-\left( 1-c_{1}\right) ^{1/\alpha }\right] .$ In that case, the
model is super-critical with $\widetilde{\phi }\left( 1\right) <1:$ the
modified branching mechanism is not regular having a probability $1-\rho
_{e} $ to take the value $\infty $ and the mean value of the offspring
distribution is infinite.

Note $\widetilde{z}_{c}=1/\left( c_{1}c_{2}\right) <\widetilde{z}%
_{*}=1/c_{1} $ if $c_{2}>1$, where $\widetilde{z}_{c}$ and $\widetilde{z}%
_{*} $ are the convergence radii of $\widetilde{\Phi }\left( z\right) $ and $%
\widetilde{\phi }\left( z\right) ,$ respectively.

For rescaled Sibuya distributions, the occupancy distribution (\ref{two})
still holds for a population of $k$ clusters in a size-$n$ forest (the
scaling factors cancel out at the denominator and numerator). However, the
law of the number $K_{n}$ of clusters is given by $\widetilde{\Phi }\left(
z\right) =c_{1}^{-1}\Phi \left( c_{1}c_{2}z\right) $%
\begin{equation*}
\widetilde{\mathbf{P}}\left( K_{n}=k\right) =\frac{n!\left[ z^{n}\right] 
\widetilde{\Phi }\left( z\right) ^{k}}{n!\left[ z^{n}\right] \left( 1-%
\widetilde{\Phi }\left( z\right) \right) ^{-1}}=\frac{c_{1}^{-k}\mathcal{C}%
_{n,k}}{\sum_{l=1}^{n}c_{1}^{-l}\mathcal{C}_{n,l}}.
\end{equation*}
where $\sum_{l=1}^{n}c_{1}^{-l}\mathcal{C}_{n,l}=n!\left[ z^{n}\right]
\left( 1-c_{1}^{-1}\left( 1-z\right) ^{\alpha }\right) ^{-1}.$ In terms of
tilted pgf's: 
\begin{equation*}
\widetilde{\mathbf{E}}\left( z^{K_{n}}\right) =\frac{s_{n}\left(
c_{1}z\right) }{s_{n}\left( c_{1}\right) }
\end{equation*}
where $s_{n}\left( z\right) :=\sum_{k=1}^{n}\mathcal{C}_{n,k}z^{k}$ are the
Sibuya degree-$n$ polynomials obeying the recursion (\ref{sibpoly}), see
Appendix \textbf{A1}.\newline

\textsl{Remarks:}

$\left( i\right) $ The regular pgf $\Phi _{\rho }\left( z\right) $ in (\ref
{f19}) of the $1$-dimensional marginal size of the forest progenies in the
thermodynamic limit is a rescaled version of $\Phi \left( z\right) $ with 
\begin{equation*}
c_{1}=\Phi \left( z_{\rho }\right) \text{ and }c_{2}=z_{\rho }/\Phi \left(
z_{\rho }\right) =\rho /\Phi ^{\prime }\left( z_{\rho }\right) .
\end{equation*}
It is thus the pgf of the progeny of a BGW process generated by the regular
pgf 
\begin{equation*}
\phi _{\rho }\left( z\right) =\frac{z_{\rho }\phi \left( \Phi \left( z_{\rho
}\right) z\right) }{\Phi \left( z_{\rho }\right) },
\end{equation*}
with $\phi $ given in (\ref{f19}) satisfying $\phi _{\rho }\left( 1\right)
=1 $ and $\Phi _{\rho }^{\prime }\left( 1\right) =1/\left( 1-\phi _{\rho
}^{\prime }\left( 1\right) \right) =\rho $ leading to $\phi _{\rho }^{\prime
}\left( 1\right) =1-1/\rho <1$ if $\rho >1$ (a sub-critical case).

An explicit case: when $\alpha =1/2$, with $\phi \left( z\right) =\frac{1}{%
2-z},$ it can be checked that $z_{\rho }=\frac{4\rho \left( \rho -1\right) }{%
\left( 2\rho -1\right) ^{2}}$, $\Phi \left( z_{\rho }\right) =\frac{2\rho -1%
}{2\rho }$ so that $\phi _{\rho }\left( 1\right) =1$ and $\phi _{\rho
}^{\prime }\left( 1\right) =z_{\rho }\phi ^{\prime }\left( \Phi \left(
z_{\rho }\right) \right) =1-1/\rho $ fulfilling the assertions.

$\left( ii\right) $ the rv with regular pgf $\widetilde{\Phi }\left(
z\right) =1-\lambda \left( 1-z\right) ^{\alpha }$, $\lambda \in \left(
0,1\right) $ is one obtained after another scaling of the Sibuya pgf. $%
\widetilde{\Phi }\left( z\right) $ cannot be the pgf of a progeny for rooted
BGW processes because $\widetilde{\Phi }\left( 0\right) \neq 0$.

\section{Sibuya rv as the progeny of an increasing tree}

A size-$n$ rooted and increasing (else recursive) labeled tree has vertices
with indices or labels $\left\{ 1,\cdots ,n\right\} $ increasing for any
path from the root to its leaves. The combinatorial version of such trees
were studied by Bergeron et al. (\textsl{1992}) where it was shown that
their counting generating function no longer obeys a functional equation,
rather a differential equation. Increasing trees are recursive ones
appearing in phylogeny, with labels of the nodes encoding their order of
appearance in the tree, and thus the chronology of evolution: the leaves of
a size$-n$ such tree are the species that can mutate to another species when
adding a new `atom' $n\rightarrow n+1$, whereas the internal nodes consist
of the species that can produce a new species in the process. A forest of
such trees consist of the different population genera. The random version of
combinatorial increasing trees appears in Drmota (\textsl{2009}), p. 15.%
\newline

Let $\Phi \left( z\right) $ solve the ordinary differential equation $\Phi
^{\prime }\left( z\right) =\varphi \left( \Phi \left( z\right) \right) $
with $\Phi \left( 0\right) =0.$ We ask if there is a generating function $%
\varphi $ with positive coefficients such that $\Phi \left( z\right)
=1-\left( 1-z\right) ^{\alpha }$, $\alpha \in \left( 0,1\right) $ (Sibuya)$.$
We find 
\begin{equation}
\varphi \left( z\right) =\alpha \left( 1-z\right) ^{-\left( 1-\alpha \right)
/\alpha },  \label{f27}
\end{equation}
with $\left[ z^{n}\right] \varphi \left( z\right) =\alpha \left[ \theta
\right] _{n}/n!>0$, $\theta :=\left( 1-\alpha \right) /\alpha >0$ and $%
\varphi \left( 0\right) =\alpha >0.$ Note that $\varphi \left( 1\right)
=\infty $ so that the generating function $\varphi $ is the one of a
positive discrete measure on the integers, now with infinite total mass.
Note also that for $\varphi \left( z\right) $ to have positive coefficients,
there is no restriction on $\alpha $. We conclude that the Sibuya pgf is the
one of the total progeny of a random increasing tree for the whole
admissible range of $\alpha $, including $\left( 0,1/2\right) $.

We now have $z=\int_{0}^{\Phi \left( z\right) }\frac{dz^{\prime }}{\varphi
\left( z\right) }$ so that $\Phi \left( z\right) $ is the inverse function
of the primitive $P\left( z\right) :=\int_{0}^{z}dz^{\prime }/\varphi \left(
z^{\prime }\right) :$ $P\left( \Phi \left( z\right) \right) =z.$ Note the
convergence radius of $\Phi \left( z\right) $ is 
\begin{equation}
z_{c}=\int_{0}^{z_{*}}\frac{dz^{\prime }}{\varphi \left( z^{\prime }\right) }%
,  \label{f33}
\end{equation}
where $z_{*}=\inf \left\{ z>0:\varphi \left( z\right) =\infty \right\} $ is
the convergence radius of $\varphi \left( z\right) $, with here, in the
Sibuya context, $z_{c}=z_{*}=1.$

The Lagrange inversion formula version for increasing trees thus states that
for all $n\geq 1$%
\begin{equation}
\left[ z^{n}\right] h\left( \Phi \left( z\right) \right) =\frac{1}{n}\left[
z^{-1}\right] h^{\prime }\left( z\right) P\left( z\right) ^{-n}.  \label{Lag}
\end{equation}
Because $P\left( 0\right) =0$, with $R\left( z\right) =z/P\left( z\right) $
(obeying $R\left( 0\right) =1/\varphi \left( 0\right) >0$), this is also 
\begin{equation*}
\left[ z^{n}\right] h\left( \Phi \left( z\right) \right) =\frac{1}{n}\left[
z^{n-1}\right] h^{\prime }\left( z\right) R\left( z\right) ^{n},
\end{equation*}
In our case, we find $P\left( z\right) =1-\left( 1-z\right) ^{1/\alpha }$
and so $R\left( z\right) =z/\left[ 1-\left( 1-z\right) ^{1/\alpha }\right] $
coincides with the offspring pgf $\phi \left( z\right) $ in (\ref{FF0}), in
the BGW setting. When $h\left( z\right) =z^{k},$ we thus get 
\begin{equation*}
\left[ z^{n}\right] \Phi \left( z\right) ^{k}=\frac{k}{n}\left[
z^{n-k}\right] R\left( z\right) ^{n}=\frac{\mathcal{C}_{n,k}}{n!}:=\frac{k!}{%
n!}\alpha ^{k}\mathcal{S}_{n,k},
\end{equation*}
which is the extension of (\ref{FF4}) to $\alpha \in \left( 0,1\right) $. We
conclude that (\ref{two}) and (\ref{FF5}) hold in the increasing trees
forests context with $\alpha $ now in the full range $\left( 0,1\right) $.
In the one-parameter approach (\ref{FF5}) with $\alpha \in \left( 0,1\right) 
$, given $n$ species are split into $k$ genera, a new incoming genus appears
with probability $\frac{\left( k+1\right) \alpha }{\alpha +n}.$ With
complementary probability, $\frac{n-k\alpha }{\alpha +n},$ there is a
connection to some node of some tree in the current size$-n$ forest, each
with size $N_{n,1}\left( k\right) $. If the connection is at a leaf of this
tree, there is no creation of a new species, rather a mutation event of the
chosen leaf to a new subsequent leaf: the leaves of the forest are thus the
currently living species that can possibly mutate to a new species. If the
connection occurs at an internal node of this tree, a new species comes to
birth: the internal nodes are the ones that can generate a new species.

In more details: an internal node of a tree is one with at least one
descendant. A leaf of a tree is one node with no descendant but with one
ascendant. From this perspective, the root of a size-$1$ new tree is a
special node, being neither an internal node nor a leaf. The internal nodes
of a tree are those species that disappeared in the past in the process of
creating new species while the leaves of a tree are those species currently
alive. Given a size-$n$ population with $k$ trees (genera) of positive sizes 
$n_{1},\cdots ,n_{k},$ for which no new genus (as a root) appears when $%
n\rightarrow n+1$, the new incoming atom connects to a size-$n_{l}$ tree
with probability $n_{l}/n,$ creating a (no) new leaf or species if it
connects to an internal node or a root (a leaf) of this size-$n_{l}$ tree.
Note that if $n_{l}=1,$ a new species appears with probability $1$. The
distribution of the number of leaves in a size$-n$ recursive increasing
random tree is addressed in Drmota (\textsl{2009}), p. 254-255.\newline

The above recursive dynamics constitutes an alternative construction of the
Yule branching speciation process [Yule (\textsl{1925})] involving Simon
distribution instead [Simon (\textsl{1955})]. See Huillet (\textsl{2022}),
Section 2.1 for a recent review where the recursive aspects of this problem
is highlighted.\newline

\textsl{Remark.} Considering the rescaled version of $\Phi \left( z\right) $%
, namely $\widetilde{\Phi }\left( z\right) =c_{1}^{-1}\Phi \left(
c_{1}c_{2}z\right) ,$ now solving $\widetilde{\Phi }^{\prime }\left(
z\right) =\widetilde{\varphi }\left( \widetilde{\Phi }\left( z\right)
\right) $ we get that 
\begin{equation*}
\widetilde{\varphi }\left( z\right) =c_{2}\varphi \left( c_{1}z\right)
=\alpha c_{2}\left( 1-c_{1}z\right) ^{-\left( 1-\alpha \right) /\alpha }
\end{equation*}
is the corresponding modified `branching mechanism'. When $c_{2}=\alpha
^{-1}\left( 1-c_{1}\right) ^{\left( 1-\alpha \right) /\alpha }$, $\widetilde{%
\varphi }\left( z\right) $ is the regular pgf of a negative-binomial rv with 
$\widetilde{\varphi }\left( 1\right) =1$. In that case, 
\begin{equation*}
\widetilde{\Phi }\left( z\right) =c_{1}^{-1}\left( 1-\left(
c_{1}c_{2}z\right) ^{\alpha }\right)
\end{equation*}
with $\widetilde{\Phi }\left( 1\right) =c_{1}^{-1}\left( 1-\left(
c_{1}c_{2}\right) ^{\alpha }\right) \leq 1$ if and only if $c_{1}\geq \alpha
/\left( 1+\alpha \right) .$ When $c_{1}=\alpha /\left( 1+\alpha \right) $,
both $\widetilde{\varphi }\left( z\right) $ and $\widetilde{\Phi }\left(
z\right) $ are regular pgf's. When $1>c_{1}>\alpha /\left( 1+\alpha \right) $%
, $\widetilde{\Phi }\left( 1\right) <1$ corresponding to a situation for
which extinction occurs with probability $c_{1}$.

\section{Appendix}

\textbf{A1.} A short reminder on generalized Stirling numbers.

Let $\left( \alpha _{1},\alpha _{2},w_{2}\right) $ be three real numbers
different from $\left( 0,0,0\right) $. Let $\left[ z:\alpha \right]
_{n}:=z\left( z+\alpha \right) ...\left( z+\left( n-1\right) \alpha \right)
=\alpha ^{n}\left[ z/\alpha \right] _{n}=\alpha ^{n}\frac{\Gamma \left(
z/\alpha +n\right) }{\Gamma \left( z/\alpha \right) },$ where $\left[
z\right] _{n}:=\left[ z:1\right] _{n}=\frac{\Gamma \left( z+n\right) }{%
\Gamma \left( z\right) }$ (the rising factorials of $z$). The generalized
Stirling numbers 
\begin{eqnarray*}
\mathcal{S}_{n,k} &\equiv &S_{n,k}\left( -\alpha _{2},-\alpha
_{1},w_{2}\right) \text{, }k=0,...,n \\
\mathcal{S}_{n,k} &\equiv &0\text{ if }k>n,
\end{eqnarray*}
are defined by the identity [Hsu and Shiue (\textsl{1998})]

\begin{equation*}
\left[ w:\alpha _{2}\right] _{n}=\sum_{k=0}^{n}\mathcal{S}_{n,k}\left[
w-w_{2}:\alpha _{1}\right] _{k}.
\end{equation*}
If $w:=w_{1}+w_{2}$, this is alternatively 
\begin{equation*}
\left[ w:\alpha _{2}\right] _{n}=\sum_{k=0}^{n}\mathcal{S}_{n,k}\left[
w_{1}:\alpha _{1}\right] _{k}.
\end{equation*}
The $\mathcal{S}_{n,k}$, $k=0,\cdots n$, obey the (Stirling's triangle)
recurrence 
\begin{equation}
\mathcal{S}_{n+1,k}=\mathcal{S}_{n,k-1}+\left( n\alpha _{2}-k\alpha
_{1}+w_{2}\right) \mathcal{S}_{n,k},  \label{rec21}
\end{equation}
with boundary conditions $\mathcal{S}_{n,0}=\left[ w_{2}:\alpha _{2}\right]
_{n}$ and $\mathcal{S}_{n,n}=1$ for all $n\geq 0$. Under certain conditions,
these numbers are positive. These numbers include the classical Stirling
numbers (signed or unsigned) of the first and second kind, when $\left(
\alpha _{1},\alpha _{2}\right) $ are specific integers, justifying the name
generalized Stirling numbers, see Huillet (\textsl{2024}).

The case $w_{2}>0$, with a combinatorial interpretation of $w_{2}$, was
discussed in M\"{o}hle (\textsl{2021}) and Huillet and M\"{o}hle (\textsl{%
2024}). We shall here assume the non-centrality parameter $w_{2}=0$ and let
(abusively) $\mathcal{S}_{n,k}:=S_{n,k}\left( -\alpha _{2},-\alpha
_{1},0\right) ,$ obtained from the general $\mathcal{S}_{n,k}$ when $%
w_{2}=0. $ These therefore obey: 
\begin{equation}
\mathcal{S}_{n+1,k}=\mathcal{S}_{n,k-1}+\left( n\alpha _{2}-k\alpha
_{1}\right) \mathcal{S}_{n,k}.  \label{rec3}
\end{equation}
These generalized Stirling polynomials $\sigma _{n}\left( u\right)
:=\sum_{k=0}^{n}\mathcal{S}_{n,k}u^{k}$ thus obey the
difference-differential recursion: 
\begin{equation*}
\sigma _{n+1}\left( u\right) =\left( u+n\alpha _{2}\right) \sigma _{n}\left(
u\right) -u\alpha _{1}\sigma _{n}^{\prime }\left( u\right) ,
\end{equation*}
making them easy to generate on a lap-top.

In fact, in this note, only the $\mathcal{S}_{n,k}$ with $\alpha
_{2}=1,\alpha _{1}=\alpha \in \left[ 0,1\right] $ will be of concern,
therefore obeying 
\begin{equation*}
\mathcal{S}_{n+1,k}=\mathcal{S}_{n,k-1}+\left( n-k\alpha \right) \mathcal{S}%
_{n,k},
\end{equation*}
with $\mathcal{S}_{n,0}=\delta _{n,0}$ and $\mathcal{S}_{n,k}>0$ if $%
k=1,\cdots ,n$ ($\mathcal{S}_{n,n}=1$).

The Dobi\'{n}ski-type formula, see Hsu and Shiue (\textsl{1998}) yields the
alternate sum expression: 
\begin{equation*}
\mathcal{S}_{n,k}:=S_{n,k}(-1,-\alpha ;0)=\frac{1}{k!\alpha ^{k}}%
\sum_{l=0}^{k}(-1)^{l}\binom{k}{l}[-l\alpha ]_{n},\text{ }k\in \{0,\ldots
,n\}.
\end{equation*}
It turns out [see Hsu and Shiue (\textsl{1998})] that, with 
\begin{equation*}
E\left( z,u\right) :=e^{u\Psi \left( z\right) },
\end{equation*}
where 
\begin{equation}
\Psi \left( z\right) =\frac{1}{\alpha }\left( 1-\left( 1-z\right) ^{\alpha
}\right) ,  \label{phi}
\end{equation}
the coefficients 
\begin{equation*}
\mathcal{S}_{n,k}:=\frac{n!}{k!}\left[ z^{n}\right] \Psi \left( z\right)
^{k},
\end{equation*}
constitute a special class of exponential Bell numbers. Hence, defining $%
\Phi \left( z\right) =\alpha \Psi \left( z\right) =1-\left( 1-z\right)
^{\alpha },$ $\mathcal{C}_{n,k}:=n!\left[ z^{n}\right] \Phi \left( z\right)
^{k}=k!\alpha ^{k}\mathcal{S}_{n,k},$ obey 
\begin{equation*}
\mathcal{C}_{n+1,k}=k\alpha \mathcal{C}_{n,k-1}+\left( n-k\alpha \right) 
\mathcal{C}_{n,k}.
\end{equation*}
Note $\mathcal{C}_{n,k}:=n!\left[ u^{k}z^{n}\right] K\left( z,u\right) $
where $K\left( z,u\right) =1/\left( 1-u\Phi \left( z\right) \right) ,$ as
explained in the main text. Defining the Sibuya polynomials by: $s_{n}\left(
u\right) =n!\left[ z^{n}\right] K\left( z,u\right) =\sum_{k=1}^{n}\mathcal{C}%
_{n,k}u^{k}$, the above recurrence translates into the
difference-differential recursion ($s_{0}\left( u\right) =1$) 
\begin{equation}
s_{n+1}\left( u\right) =\left( n+\alpha u\right) s_{n}\left( u\right)
+\alpha u\left( u-1\right) s_{n}^{\prime }\left( u\right) \text{, }n\geq 0%
\text{,}  \label{sibpoly}
\end{equation}
with $s_{n+1}\left( 1\right) =\left( n+\alpha \right) s_{n}\left( 1\right) ,$
hence $s_{n}\left( 1\right) =\left[ \alpha \right] _{n}.$ It can be checked,
for example on $s_{3}\left( u\right) $, that these polynomials are not in
the class of polynomials with real zeros, see Pitman (\textsl{1997}).\newline

\textbf{A2.} We herewith work out the connection of occupancies for forests
of critical Sibuya trees with the Ewens-Pitman sampling formula.

With $\alpha \in \left[ 0,1\right] $ and $\theta >-\alpha $, let $S_{n}$ be
a pure-birth Markov chain obeying $S_{0}=0$ and for $n\geq 0,$

\begin{equation*}
\mathbf{P}(S_{n+1}=k)=\frac{\theta +\left( k-1\right) \alpha }{\theta +n}%
\mathbf{P}(S_{n}=k\,|\,S_{n}=k-1)+\frac{n-k\alpha }{\theta +n}\mathbf{P}%
(S_{n+1}=k\,|\,S_{n}=k).
\end{equation*}
Equivalently, 
\begin{equation*}
\mathbf{P}(S_{n+1}=k+1\,|\,S_{n}=k)=1-\mathbf{P}(S_{n+1}=k\,|\,S_{n}=k):=%
\frac{\theta +k\alpha }{\theta +n},\text{ }n\in \Bbb{N}_{0}.
\end{equation*}
In this case $S_{n}$ coincides with the number of occupied tables in the
(Chinese) restaurant process [Pitman (\textsl{2006})] after $n$ customers
have entered the restaurant. The transition probabilities of the random walk 
$(S_{n},n\in \Bbb{N}_{0})$ not only depend on the `time' $n$ but also on the
current state $S_{n}=k$. The integrated distribution of $S_{n}$ can be
expressed as [see Pitman (\textsl{2006}), (3.11) p. 66] 
\begin{equation*}
\mathbf{P}(S_{n}=k)=\frac{[\theta :\alpha ]_{k}}{[\theta ]_{n}}%
S_{n,k}(-1,-\alpha ;0),\text{ }k\in \{0,\ldots ,n\},
\end{equation*}
where $[\theta :\alpha ]_{0}:=1$, $[\theta :\alpha
]_{k}:=\prod_{l=0}^{k-1}(\theta +l\alpha )$ for $k\in \Bbb{N}$.

$S_{n}$ counts the number of distinct species in a sample of size $n$ taken
from Pitman and Yor (\textsl{1997}) two-parameter stick-breaking $\mathrm{PD}
$($\alpha ,\theta $)-partition of the unit interval, extending the Ewens
case (for which $\alpha =0$). When $\alpha >0$ ($\alpha =0$), species
abundances are ranked in decreasing order, with small or rare abundances
decaying algebraically (geometrically) fast. We refer the reader to Chapter
3 of Pitman's lecture notes [Pitman (\textsl{2006})] for further information
on the two-parameter model and to Yamato and Sibuya (\textsl{2000}) and
Yamato, Sibuya and Nomachi (\textsl{2001}) for some further related works on
sampling.

For $\theta >0$ and $\alpha \in \left( 0,1\right) $, as $n\to \infty $, $%
n^{-\alpha }S_{n}$ converges almost surely and in $L^{q}$ for any $q>0$ to a
limiting random variable $W_{\alpha ,\theta }$ being two-parameter
Mittag--Leffler distributed, for which 
\begin{equation*}
\mathbf{E}\left( W_{\alpha ,\theta }^{q}\right) =\frac{\Gamma \left( \theta
\right) }{\Gamma \left( \theta +q\alpha \right) }\frac{\Gamma \left( \theta
/\alpha +q\right) }{\Gamma \left( \theta /\alpha \right) },\text{ }q>-\theta
/\alpha .
\end{equation*}
See Pitman (\textsl{2006}), Theorem 3.8. $W_{\alpha ,\theta }$ captures the
biodiversity of the partition.

In this sampling context, the following occupancy distribution arising from
Ewens--Pitman sampling formula holds: 
\begin{eqnarray*}
&&\hspace{-15mm}\mathbf{P}(N_{n}(1)=n_{1},\ldots ,N_{n}(k)=n_{k},S_{n}=k) \\
&=&\frac{n!}{k!}\frac{[\theta :\alpha ]_{k}}{[\theta ]_{n}}\prod_{l=1}^{k}%
\frac{[1-\alpha ]_{n_{l}-1}}{n_{l}!} \\
&=&\binom{n}{n_{1}\cdots n_{k}}\frac{[\theta :\alpha ]_{k}}{k![\theta ]_{n}}%
\prod_{l=1}^{k}[1-\alpha ]_{n_{l}-1}
\end{eqnarray*}
where to the right-hand side we recognize the joint distribution of the box
occupancies $N_{n}(l)^{\prime }$s when uniform sampling is from $\mathrm{PD}%
(\alpha ,\theta )$, observed in an arbitrary way (independently of the
sampling mechanism) and with total sample size $\sum_{l=1}^{k}n_{l}=n.$
Consequently, 
\begin{eqnarray*}
\mathbf{P}(N_{n,k}(1) &=&n_{1},\ldots ,N_{n,k}(k)=n_{k}):=\mathbf{P}%
(N_{n}(1)=n_{1},\ldots ,N_{n}(k)=n_{k}\mid S_{n}=k) \\
&=&\binom{n}{n_{1}\cdots n_{k}}\frac{\alpha ^{k}}{k!\alpha ^{k}\mathcal{S}%
_{n,k}}\prod_{l=1}^{k}[1-\alpha ]_{n_{l}-1},
\end{eqnarray*}
independent of $\theta $. Recalling $\mathcal{C}_{n,k}=k!\alpha ^{k}\mathcal{%
S}_{n,k}$, this expression coincides with (\ref{two}).

Finally, the result (\ref{FF6}) follows from the observation that, when $%
\theta =\alpha $, the dynamics of $K_{n}$ and $S_{n}$ coincide: the
recursive formation of critical Sibuya forests matches with the discovery of
new species in the sampling process from $\mathrm{PD}$($\alpha ,\alpha $).%
\newline

\textbf{Acknowledgments:} This work benefited from the support of the Chair
``Mod\'{e}lisation Math\'{e}matique et Biodiversit\'{e}'' of Veolia-Ecole
Polytechnique-MNHN-Fondation X. The author also warmly thanks Professor
Martin M\"{o}hle, from the University of T\"{u}bingen, for his careful
reading of the manuscript and his kind advice to improve an early version of
it.

\end{document}